\documentclass[12pt]{article}
\usepackage{graphicx,amscd,amssymb}
\begin{document}

\newcommand{\PSbox}[3]{\mbox{\rule{0in}{#3}\includegraphics{#1}\hspace{#2}}}
\newtheorem{theorem}{Theorem}[section]
\newtheorem{definition}[theorem]{Definition}
\newtheorem{example}[theorem]{Example}
\newtheorem{lemma}[theorem]{Lemma}
\newtheorem{proposition}[theorem]{Proposition}
\newtheorem{corollary}[theorem]{Corollary}
\newtheorem{remark}[theorem]{Remark}
\newtheorem{conjecture}[theorem]{Conjecture}

\newcommand{\A}{{\cal A}}
\newcommand{\B}{{\cal B}}
\newcommand{\C}{{\cal C}}
\newcommand{\E}{{\cal E}}
\newcommand{\F}{{\cal F}}
\newcommand{\Ho}{{\cal H}}
\newcommand{\M}{{\cal M}}
\newcommand{\N}{{\cal N}}
\newcommand{\V}{{\cal V}}
\newcommand{\Lo}{{\cal L}}
\newcommand{\X}{{\cal X}}
\newcommand{\al}{\alpha}
\newcommand{\Ann}{\mbox{\rm Ann}}

\newcommand{\cn}{{\bf {\rm C}}
\hspace{-.4em}      {\vrule height1.5ex width.08em depth-.04ex}
\hspace{.3em}}
\newcommand{\Hom}{\mbox{\rm Hom}}
\newcommand{\Der}{\mbox{\rm Der}(S)}
\newcommand{\Char}{\mbox{\rm Char}}
\newcommand{\gl}{{g_{\lambda}}}
\newcommand{\Gr}{\mbox{\rm Gr}}
\newcommand{\ints}{{\sf Z}\hspace{-.36em}{\sf Z}}
\newcommand{\la}{\lambda}
\newcommand{\om}{\omega}
\newcommand{\Om}{\Omega}
\newcommand{\Op}{\Omega^{p}}
\newcommand{\OpA}{\Omega^{p}(\A)}

\newcommand{\p}{\partial}
\newcommand{\PD}{\mbox{\rm pd\,}}
\newcommand{\Pl}{{\Phi_{\lambda}}}
\newcommand{\pl}{{\phi_{\lambda}}}
\newcommand{\Poin}{\mbox{\rm Poin}}
\newcommand{\PoinO}{\Poin(\Omega^{*}(c\A); u, }
\newcommand{\PoinGr}{\Poin(\Gr\Omega^{*}(\A); u, }
\newcommand{\rk}{\mbox{\rm rank}}
\newcommand{\proof}{{\bf Proof.~}}
\newcommand{\qed}{~~\mbox{$\Box$}}

\newcommand{\ra}{{\rightarrow}}
\newcommand{\rn}{{\rm I}\hspace{-.2em}{\rm R}}

\newcommand{\scn}{\scriptsize\cn}

\newcommand{\we}{\wedge}
\newcommand{\ar}{\buildrel d_{\la}\over\rightarrow}

\title{  Hodge decomposition of Alexander invariants.}
\author  {
{\sc Anatoly Libgober}\\
{\small\it Department of Mathematics, University of Illinois,
Chicago, Ill 60607}}

%\date{to be included}

\maketitle
\begin{abstract}
{Multivariable Alexander invariants of algebraic links 
calculated in terms of algebro-geometric invariants 
links (polytopes and ideals of quasiadjunction). 
The relations with log-canonical divisors, the multiplier ideals
and a semicontinuity property of polytopes of quasiadjunction
is discussed.} 
\end{abstract}

\section{Introduction}
\bigskip

This paper is concerned with two interrelated issues. The first is an algebraic
calculation of multivariable 
Alexander invariants of the fundamental groups of links of plane curves 
singularities $f_1(x,y) \cdot \cdot \cdot f_r(x,y)$ 
with several branches or more precisely the corresponding 
{\it characteristic varieties} and 
the effect of Hodge theory on the structure of these invariants.

Characteristic varieties are attached to a topological space 
$X$ which  fundamental 
group admits a surjection on a free abelian group ${\bf Z}^r$.
They can be defined as follows (cf. \cite{topandappl}). 
To  a surjection onto ${\bf Z}^r$ corresponds the abelian 
cover $\tilde X$ for which ${\bf Z}^r$ is the group of 
deck transformations.  Consider 
$H_1(\tilde X,{\bf C})$ as a module over the group ring of 
${\bf Z}^r$ and let $\Phi$ be the map in a presentation 
${\bf C}[{\bf Z}^r]^m {\buildrel \Phi \over \longrightarrow} 
{\bf C}[{\bf Z}^r]^n \longrightarrow H_1(\tilde X,{\bf C})
\longrightarrow 0$
of ${\bf C}[{\bf Z}^r]$-module $H_1(\tilde X,{\bf C})$ by 
generators and relators. The $i$-th Fitting ideal of 
$H_1(\tilde X,{\bf C})$ is the ideal 
generated by the minors of order $n-i+1$ in the matrix
$\Phi$.
The ring ${\bf C}[{\bf Z}^r]$ can be 
viewed as the ring of regular functions on the torus 
${{\bf C}^*}^r$ and the $i$-th characteristic variety of $X$
is defined as the zero set in this torus of the $i$-th Fitting 
ideals of the module $H_1(\tilde X,{\bf C})$.

In the case when $X$ is the complement to a link $L$ 
with $r$-components the group
$H_1(S^3-L,{\bf Z})$ is isomorphic to ${\bf Z}^r$.   
The first Fitting ideal is a product 
of a power of the maximal ideal 
$\cal M$ of the point in ${{\bf C}^*}^r$
corresponding to the identity element of the torus  
and a principal ideal generated by a polynomial
$\Delta(t_1,...t_r)$ (cf. \cite{Crowell}). The latter is 
called the Alexander polynomial of $L$.
In particular, the first characteristic variety has  
codimension 1 in ${\bf C}^r$. For $i >1$ the characteristic varieties
may have higher codimension. 
For any positive $r'< r$ 
and surjection $H_1(S^3-L,{\bf Z}) \rightarrow {\bf Z}^{r'}$ 
one has the Alexander polynomial of $r'$ variables 
and characteristic varieties in ${{\bf C}^*}^{r'}$ which can 
be found from those of $r$ variables (cf. \cite{turaev}).
One variable Alexander polynomial can be obtained from this construction
applied to surjection $H_1(S^3-L,{\bf Z}) \rightarrow {\bf Z}$ 
given by evaluation the total linking number of loops 
with $L$. 

If $r=1$ {\it all} Fitting ideals are principal and their generators 
yield a sequence of polynomials $\Delta_i \in {\bf C}[t,t^{-1}]$
(defined up to a unit of the latter) 
such that $\Delta_{i+1} \vert \Delta_i$. 
In the case when $L$ is the link of a singularity $f=0$, 
both, the algebraic calculation of Alexander invariants 
and their relation to the Hodge theory are well known.
Indeed, the corresponding infinite cover of $S^3-L$ is cyclic, it 
can be identified with the Milnor fiber of the singularity,
and in this identification the deck transformation is the monodromy 
operator of the singularity. The group $H^1$ of the Milnor fiber
supports a mixed Hodge structure with weights 0,1 and 2, 
with the identification $N: W_2/W_1 \rightarrow W_0$ 
given by the logarithm of an appropriate power of 
the monodromy (cf. \cite{oslo})
All Hodge groups are invariant under the action of the 
semisimple part of the latter. 
Let $h^{p,q}_{\zeta}$ (cf. \cite{oslo}) be the dimension of the eigenspace
of this semisimple part acting on the space $H^{p,q}$.
These numbers determine the Jordan form of the monodromy as follows.
The size of the Jordan blocks 
of the monodromy does not exceed 2 
and the number of blocks corresponding to
an eigenvalue $\zeta$ of size $1 \times 1$ (resp. $2 \times 2$)
is equal to $h^{1,0}_{\zeta}+h^{0,1}_{\zeta}$ (resp. $h^{0,0}_{\zeta}$). 
As a consequence:  
$$\Delta_i=\prod_{(\zeta)}(t-\zeta)^{a_{\zeta,i}}$$
where 

$$a_{\zeta,i}=\cases { h_{\zeta}^{1,0}+h_{\zeta}^{0,1}+2h_{\zeta}^{0,0}
-2(i-1) 
&  if  $1 \le i \le h_{\zeta}^{0,0}$  \cr
h_{\zeta}^{1,0}+h_{\zeta}^{0,1}-(i-1-h_{\zeta}^{0,0}) 
&  if  $ h_{\zeta}^{0,0} < i \le h_{\zeta}^{0,0}+h_{\zeta}^{1,0}+
h_{\zeta}^{0,1}$  \cr
0 &  if $i >  h_{\zeta}^{0,0}+h_{\zeta}^{1,0}+
h_{\zeta}^{0,1}$  \cr }$$

All $\Delta_i$  can be calculated
algebraically in terms of a resolution of the singularity.
For $\Delta_1$ this follows form A'Campo's formula for the 
$\zeta$-function of the monodromy (cf. \cite{A'Campo}) and 
for $\Delta_i$ and $i \ge 2$ from Steenbrink's calculation of the Hodge 
numbers of Mixed Hodge structure on the cohomology of the 
Milnor fiber (cf. \cite{oslo}). 

Multivariable Alexander polynomials were studied extensively 
from topological point of view. They can be found either from
a presentation of the fundamental group $\pi_1(S^3-L)$ 
using Fox calculus or using an iterative procedure based 
on the fact that algebraic links are iterated torus links
(cf. \cite{sumners} \cite{EN}, cf. also \cite{Sabbah} 
where an upper bound for the 
set of zeros of the multivariable Alexander polynomial 
was obtained algebraically). 

In this paper we describe an algebraic procedure for calculating
the characteristic varieties. In fact   
we study a finer than characteristic variety invariant of singularity 
having a given 
algebraic link. This invariant is  
a collection of polytopes in ${\bf R}^r$.
These polytopes are equivalent to the local polytopes
of quasiadjunction introduced in \cite{Abcov} but are more
convenient in the local case. 
We show that the characteristic varieties are algebraic 
closures of the images of the faces of these polytopes of 
quasiadjunction under exponential map: $exp: {\bf R}^r \rightarrow
{{\bf C}^*}^r$. In the case $r=1$ polytopes of quasiadjunction are
segments having $1$ as the right end and faces of quasiadjunction are
points in [0,1] which are the left ends of these segments. 
These points are elements of the Arnold-Steenbrink spectrum of the 
singularity (the analogy is going further: we prove 
in \ref{semicontinuity} a semicontinuity 
property of faces of quasiadjunction 
extending one of well known semicontinuity properties of spectra
(\cite{semicontinuity}, \cite{varchenko}).
Description of characteristic varieties via faces of quasiadjunction
is obtained by 
expressing the mixed Hodge structure on the cohomology of the 
abelian covers of $S^3$ branched over the link of singularity   
in terms of certain ideals in the local ring of singularity
(ideals of quasiadjunction). These ideals are the key ingredient
in the description of the characteristic varieties of the 
fundamental groups of the complements to curves 
in ${\bf P}^2$ (cf. \cite{wall}, \cite{Abcov} ) and are 
generalizations of the 
ideals of quasiadjunction in $r=1$ case (cf. \cite{arcata1}, \cite{loser}). 
The ideals of quasiadjunction (in both $r=1$ and $r >1$) cases 
are closely related to more recently introduced 
multiplier ideals (cf. \cite{nadel} and remark \ref{multiplier}). 
On the other hand, following an idea of J.Kollar (cf. \cite{kollar}), 
we show how log-canonical thresholds of certain divisors can be found 
in terms of polytopes of quasiadjunction studied here and in \cite{Abcov}
(cf. section (\ref{thresh})).

In  section \ref{properties} we point out the effect of the Hodge theory 
on the homology of the infinite abelian covers of the complements 
to links. Note that these homology groups typically are infinite dimensional. 
We show that the intersections of the torus of unitary characters 
of ${\bf Z}^r$ with irreducible components of the space of characters 
appearing in the representation of ${\bf Z}^r$ 
on $H_1(\tilde X,{\bf C})$ have natural decomposition into a union of two connected subsets
compatible with the Hodge decomposition of the finite abelian covers
and illustrate such decomposition by explicit examples. I thank
J.Cogolludo for sharing with me his results on Fox calculus 
calculations in example 2. 
\par Finally note, that polytopes and ideals of 
quasiadjunction considered here have a natural generalization for 
arbitrary hypersurface singularities with applications generalizing 
\cite{homotopy}. We shall return to this elsewhere.

\section{Invariants of singularities.}

\subsection{Characteristic varieties of algebraic links}
\label{translatedtori}

Let $X$, as in Introduction  (cf. also \cite{topandappl}, \cite{Abcov},
 \cite{Sabbah}), be a finite CW complex such that $H_1(X,{\bf Z})={\bf Z}^r$.
Let $t_1,..,t_r$ be a system of generators of the latter.
The homology $H_1(\tilde X,{\bf C})$
of the universal abelian cover $\tilde X$
has a structure of 
$\Lambda={\bf C}[H_1(X,{\bf Z})]=
{\bf C}[t_1,t_1^{-1},..,t_r,t_r^{-1}]$-module. 
Let $F_i(X)$ be the $i$-th Fitting ideal of $H_1(\tilde X,{\bf C})$
considered as a $\Lambda$-module i.e. the ideal generated by 
$(n-i+1) \times (n-i+1)$ minors of the matrix of a map $\Phi: 
\Lambda^m \rightarrow \Lambda^n$ such that ${\rm Coker} \Phi=
H_1(\tilde X,{\bf C})$. 

We shall view  ${\bf C}[t_1,t_1^{-1},..,t_r,t_r^{-1}]$ 
as a ring of regular functions on a torus ${{\bf C}^*}^r$ so that 
the set of zeros of an ideal $F_i(X)$ is a subvariety of  
${{\bf C}^*}^r$
denoted $V_i(X)$. The maximal possible $i$ is called the {\it depth} of $V_i$.
A translated subgroup of ${{\bf C}^*}^r$ is an irreducible component 
of an intersection of codimension one submanifolds given by 
$t_1^{l_1} \cdot \cdot \cdot t_r^{l_r}=\lambda (l_i \in {\bf Z})$.
The following is a local analog of results of \cite{arapura}.

\begin{proposition} The characteristic varieties of 
algebraic links are unions of translated subgroups.
\end{proposition}

\noindent \proof Recall that an algebraic link can be obtained 
from a trivial knot by iteration
of cablings $X_1 \rightarrow ... \rightarrow X_N$ 
where $X_1$ is a complement to a small tube about the unknot and the
complement $X_n$ to a link $L_n$  
is obtained from the complement $X_{n-1}$ to the link $L_{n-1}$ 
via replacing $X_{n-1}$ 
by a space $X_{n-1} \cup_{\partial T_n} Y_{n-1}$. Here $Y_{n-1}$
is one of two standard model spaces: complement in a torus $T_n$
either to a torus knot or to the union of the axis of the torus and 
the torus 
knot (cf. \cite{sumners}). The union is taken 
by identifying the boundary of a tube about $L_{n-1}$ with 
$\partial T_n$.
It follows from \cite{sumners} that for the homology of 
the universal abelian cover $\tilde X_n$ we have:
 $$H_1(\tilde X_n)=H_1(\tilde X_{n-1})\oplus H_1(\tilde Y_{n-1})$$
(cf. p.118 and p.119 in \cite{sumners} for two possible cases of cablings).  
We can assume by induction that the characteristic varieties 
of $\tilde X_{n-1}$ and $\tilde Y_{n-1}$ are the unions of translated
subgroups. We have $V_k(X_n)=Supp(\Lambda^k H_1(\tilde X_n))$ 
(cf. {\cite{Abcov}). $\Lambda^k H_1(\tilde X_n)$ has a filtration 
with successive factors $\Lambda^i H_1(\tilde X_{n-1}) \otimes 
\Lambda^{k-i} H_1(\tilde Y_{n-1})$. Hence, by \cite{serre},
$V_k(\tilde X_n)=\cap_i V_i(\tilde X_{n-1}) \cup V_{k-i}(\tilde
Y_{n-1})$. In particular, if $V_i(\tilde X_{n-1})$ and $ V_{k-i}(\tilde
Y_{n-1})$ are unions of translated subgroups then so is $V_k(\tilde X_n)$.

\subsection{Ideals of log-quasiadjunction}

Let $B$ be a small ball about the origin $O$ in ${\bf C}^2$ and let 
$C$ be a germ of a plane curve having at $O$  
singularity with $r$ branches.  
Let $f_1(x,y) \cdot \cdot \cdot f_r(x,y)=0$ 
where $f_i=0$ is irreducible be a local equation of this curve.
 An abelian cover of type
$(m_1,...,m_r)$ of $\partial B$ (resp. $B$)
is the branched cover of $\partial B$ (resp. $B$) 
corresponding to a 
homomorphism $\pi_1(\partial B -\partial B \cap C) \rightarrow 
{\bf Z}/m_1{\bf Z} \oplus ... \oplus {\bf Z}/m_r{\bf Z}$ (resp. the cone 
over the abelian cover of $\partial B$). 
Such cover of $\partial B$ 
is the link of complete intersection surface singularity:
\begin{equation}V_{m_1,..,m_r}: 
\ \ \ z_1^{m_1}=f_1(x,y),..,z_r^{m_r}=f_r(x,y)
\label{completeintersection}
\end{equation} 
The covering map is given by $p: (z_1,...,z_r,x,y) \rightarrow (x,y)$.

\indent {\it An ideal of quasiadjunction} 
of type $(j_1,..,j_r \vert m_1,..,m_r)$ (cf. \cite{wall} \cite{Abcov})
is the ideal in the local ring of the singularity of
$C$ (i.e. $O \in {\bf C}^2$) consisting of germs $\phi$ such that 
the 2-form:
\begin{equation}
\omega_{\phi}={{\phi z_1^{j_1} \cdot \cdot \cdot z_r^{j_r} dx  \wedge 
dy} \over  {z_1^{m_1-1} \cdot \cdot \cdot z_r^{m_r-1}}}
\label{form}
\end{equation}
extends to a holomorphic form on a resolution of the singularity 
of the abelian cover of a ball $B$ of type $(m_1, ...., m_r)$
i.e. a resolution of (\ref{completeintersection})
(we suppress dependence of $\omega_{\phi}$ on $j_1,..j_r, m_1,...m_r$).
In other words, $\phi z_1^{j_1} \cdot \cdot \cdot z_r^{j_r}$ belongs to 
the adjoint ideal of (\ref{completeintersection}) (cf. \cite{merle}).
In particular the condition on $\phi$ is independent of resolution.
We always shall assume that 
$0 \le j_1 < m_1,..,0 \le j_r < m_r$.

\smallskip 
\noindent 
{\it An ideal of log-quasiadjunction} 
(resp. {\it an ideal of weight one} {\it log-quasiadjunction}) 
of type $(j_1,..,j_r \vert m_1,..,m_r)$
is the ideal in the same local ring 
consisting of germs $\phi$ such that $\omega_{\phi}$ 
extends to a log-form (resp. weight one log-form) on a resolution 
of the singularity of the same abelian cover. Recall (cf. \cite{deligne})
that a holomorphic 2-form is weight one log-form if it is a 
combination of forms having
poles of order at most 
one on each component of the exceptional divisor and  
having no poles on a pair of intersecting components. 
These ideals are also independent of a resolution. 
This follows from the following.
Let $\omega$ be a holomorphic $n$-form on a complex space 
$X \subset {\bf C}^N,{\rm dim} X=n$
with isolated singularity at $x \in X$ and $B_r$ be a ball 
about $x$ in ${\bf C}^N$ having a radius $r$. $\omega$ extends to a form 
of weight $k$ on a resolution of $X$, for which the exceptional 
divisor has at worse normal crossings, if and only if for sufficiently
 small $R \gg r >0$ one has 
$$  \int_{B_R-B_r} \omega \wedge \bar \omega   < C {\vert {\rm log} \ 
 r
\vert}^k$$
In particular for $k=0$ one obtains Lemma 1.3 (ii) from \cite{merle}.
General case follows, for example, by 
interpreting the above integral as an integral over the neighborhood of
the exceptional locus in a resolution of $x \in X$ and reducing it
to the integral of over the boundary of $B_R-B_r$. Local 
calculations near intersection of $k$ components shows 
that the contribution   $  \int_{r \le z_1 \le 1,....
r \le z_k \le 1, 0 \le z_{k+1} \le 1,  0 \le z_{n} \le 1} 
{{dz_1 ... dz_n  \wedge \bar {dz_1} ...  \bar {dz_n}} \over 
{z_1 \bar {z_1} \cdot \cdot \cdot z_n \bar {z_n}}}  < C {\vert {\rm log}\  r
\vert}^k $ which yields the estimate as above
(similarly to \cite{merle}). 
This characterization gives independence of particular resolution
both  the ideals of 
log-quasiadjunction and ideal of weight one log-quasiadjunction.

\smallskip  It is shown in \cite{Abcov} that an ideal of quasiadjunction 
${\cal A}(j_1,..,j_r \vert m_1,...,m_r)$ is determined by the vector:
\begin{equation}
({{j_1+1} \over {m_1}},....,{{j_r+1} \over
  {m_r}}).
\label{vector}
\end{equation}
 This is also the case for the ideals of 
log-quasiadjunction and weight one log-quasiadjunction. 
Indeed, these ideals can be described as follows. 
For a given embedded resolution $\pi: V \rightarrow {\bf C}^2$ of the 
germ $f_1 \cdot \cdot \cdot f_r=0$ with the exceptional curves 
$E_1,..,E_k,...,E_s$ let $a_{k,i}$ (resp. $c_k$, resp. $e_k(\phi)$)
 be the multiplicity of the pull back on $V$ of 
$f_i$ ($i=1,..,r$) (resp. $dx \wedge dy$, resp. $\phi$)
along $E_k$.
Then $\phi$ belongs to the ideal of quasiadjunction of type 
   $(j_1,..,j_r \vert m_1,..,m_r)$ if and only if for any $k$, 
\begin{equation}
 a_{k,1}{{j_1+1} \over {m_1}}+....,+a_{k,r}{{j_r+1} \over {m_r}}
 > a_{k,1}+...+a_{k_r}-e_k(\phi)-c_k-1
\label{quasiadj}
\end{equation}
(cf. \cite{Abcov}).
Similar calculation shows that 
a germ $\phi$ belongs to the ideal of log-quasiadjunction
corresponding to 
$(j_1,..,j_r \vert m_1,..,m_r)$ if and only if the 
inequality 
\begin{equation}
a_{k,1}{{j_1+1} \over {m_1}}+....,+a_{k,r}{{j_r+1} \over {m_r}}
 \ge  a_{k,1}+...+a_{k_r}-e_k(\phi)-c_k-1
\label{logquasiadj}
\end{equation}
is satisfied for any $k$.
Moreover,  a germ $\phi$ belongs 
to the ideal of weight one log-quasiadjunction if and
only if this germ is a linear combination of germs $\phi$ 
satisfying inequality (\ref{logquasiadj})
for any collection 
of $k$'s such that corresponding components do not 
intersect 
and satisfying the inequality  (\ref{quasiadj}) 
for $k$ outside of this collection. 
We shall denote the ideal of quasiadjunction (resp. log-quasiadjunction, 
resp. weight one log-quasiadjunction) corresponding to 
$(j_1,..,j_r \vert m_1,..,m_r)$ as ${\cal A} (j_1,..,j_r \vert
m_1,..,m_r)$ (resp. ${\cal A}'' (j_1,..,j_r \vert
m_1,..,m_r)$, resp. ${\cal A}' (j_1,..,j_r \vert
m_1,..,m_r)$). We have: 
$${\cal A} (j_1,..,j_r \vert
m_1,..,m_r) \subseteq  
{\cal A}' (j_1,..,j_r \vert
m_1,..,m_r) \subseteq 
{\cal A}'' (j_1,..,j_r \vert m_1,..,m_r)$$

\par \noindent
Recall that both (\ref{quasiadj}) and (\ref{logquasiadj}) follow 
from the following calculation (cf. \cite{Abcov} section 2
for complete details).
One can use the normalization of the fiber product
$\widetilde V_{m_1,...,m_r}=
V \times_{{\bf C}^2} V_{m_1,..,m_r}$
as a resolution of singularity (\ref{completeintersection})
in the category of manifolds with quotient singularities 
(cf. \cite{lipman}).
We have: 
\begin{equation}
\matrix{  \tilde V_{m_1,..,m_r} & \buildrel \tilde p \over 
                                                \rightarrow & V \cr
               \tilde \pi  \downarrow & &  \pi \downarrow \cr
               V_{m_1,..,m_r}&  \buildrel p \over 
                     \rightarrow & {\bf C}^2 \cr }
\label{resolution}
\end{equation}        
Preimage of the exceptional divisor of $V \rightarrow {\bf C}^2$
in $\widetilde V_{m_1,...,m_r}$ forms a divisor with normal crossings
(cf. \cite{oslo}), though preimage of each component is reducible 
in general (in which case irreducible components 
above each exceptional curve do not intersect
\footnote{${}^*${If the Galois group $G$ of $\tilde p$ is abelian
(as we always assume here) and, in particular, is the quotient of 
$H_1(B-C \cap B,{\bf Z})$ 
then the Galois group of ${\tilde p}^{-1}(E_i) \rightarrow E_i$ 
is $G/(\gamma_i)$ 
where for an exceptional curve $E_k$, $\gamma_k$ is the image in 
the Galois group of the homology class of the boundary of a small 
disk transversal to $E_k$ in $V$. The components of 
${\tilde p}^{-1}(E_i)$ correspond to the elements of 
$G/(\gamma_i, ...\gamma_l ..)$ where $l$ runs through indices 
of all exceptional curves intersecting $E_i$, while 
$\tilde p_i$ restricted on each component has 
$(\gamma_i, ...\gamma_l ..)/(\gamma_i)$ as the Galois group.
The points ${\tilde p}^{-1}(E_i \cap E_j)$ 
correspond to the elements of $G/(\gamma_i,\gamma_j)$
and the points of ${\tilde p}^{-1}(E_i \cap E_j)$
belonging to a fixed component correspond to cosets
in $ (\gamma_i, ...\gamma_l ..)/(\gamma_i, \gamma_j).$ } }  ).
The order of the vanishing of $\omega_{\phi}$ on 
$\widetilde V_{m_1,...,m_r}$ along $E_k$ is equal to:
\begin{equation}
\Sigma_{i=1}^{i=r}
 ( j_i-m_i+1) {{m_1 \cdot \cdot \cdot \hat m_i \cdot \cdot \cdot m_r
 \cdot a_{k,i}} \over {g_{k,1} \cdot \cdot \cdot g_{k,r}s_k}}+
{{m_1 \cdot \cdot \cdot m_r \cdot ord_{E_k}(\pi^*(\phi))} \over {g_{k,1} \cdot
\cdot \cdot g_{k,r} \cdot s_k}}+
\label{orderofzero}
\end{equation}
$$+{{c_k \cdot m_1 \cdot \cdot \cdot m_r} \over {g_{k,1} \cdot \cdot \cdot
 g_{k,r} \cdot s_k}}+ {{m_1 \cdot \cdot \cdot m_r} \over {g_{k,1} \cdot
\cdot \cdot g_{k,r} \cdot s_k}}-1 $$
where $g_{k,i}=g.c.d.(m_i,a_{k,i})$ and 
$s_k=g.c.d.(...,{{m_i} \over {g_{k,i}}},...)$.

A consequence of (\ref{orderofzero}) is that $\omega_{\phi}$ has 
order of pole equal to one (resp. zero) along the component $E_k$ of the 
above resolution if an only if for such $\phi$ one has 
equality in (\ref{logquasiadj}) (resp. (\ref{quasiadj}) is satisfied).

\begin{proposition} 1. Let ${\cal A}''$ be an ideal of log-quasiadjunction. 
There is unique polytope ${\cal P}({\cal A}'')$ such that
a vector $({{j_1+1} \over {m_1}},...
 {{j_r+1} \over {m_r}}) \in {\cal P}({\cal A}'')$ 
if and only if  ${\cal A}'' (j_1,..,j_r \vert
m_1,..,m_r)$ contains ${\cal A}''$ 
\footnote{$(^*)${ i.e. a subset in ${\bf R}^r$ given by a set 
of linear inequalities $L_s \ge k_s$. 
We say that an affine hyperplane in ${\bf R}^r$ supports a codimension one face
of a polytope if the intersection this hyperplane with the boundary 
of the polytope has dimension $r-1$. 
A face of a polytope is the intersection of a supporting face of 
the polytope with the boundary. A codimension one face of a polytope 
in ${\bf R}^r$ is a polytope of dimension $r-1$.
By induction one obtains faces of arbitrary codimension
for original polytope (for $r=3$ those are called edges and vertices). 
The boundary of the polytope is the union of its faces.}}.
\par 2. The set of vectors (\ref{vector}) for which 
${\cal A} (j_1,..,j_r \vert m_1,..,m_r) \ne$ 
\newline ${\cal A}''(j_1,..,j_r \vert m_1,..,m_r)$ is a dense 
subset in the boundary of the polytope having as its closure 
a union of faces of such a polytope.
The closure of the set of vectors (\ref{vector}) 
for which  ${\cal A}' (j_1,..,j_r \vert m_1,..,m_r) \ne 
{\cal A''}(j_1,..,j_r \vert m_1,..,m_r)$
is also a union of certain faces of such a polytope. 
\par 3. The ideal ${\cal A}(j_1,..,j_r \vert  m_1,...,m_r)$ 
 (resp. ${\cal A}'(j_1,..,j_r \vert  m_1,...,m_r)$ and 
\newline ${\cal A}''(j_1,..,j_r \vert  m_1,...,m_r)$) is independent of 
array $(j_1,..,j_r \vert m_1,..,m_r)$ as long as the vector (\ref{vector}) 
varies within the interior of the same face of quasiadjunction.

\end{proposition}

\noindent We shall call the above faces {\it the faces of quasiadjunction} 
(resp. {\it weight one faces of quasiadjunction}). ${\cal A}_{\Sigma}$ will 
denote ${\cal A}(j_1,...,j_r \vert m_1,...,m_r)$ with corresponding
vector (\ref{vector}) belonging to the interior of a face of quasiadjunction $\Sigma$
(similarly for ${\cal A}_{\Sigma}'$ and ${\cal A}_{\Sigma}''$).

\bigskip 
\noindent \proof First let us describe 
the inequalities defining the polytope
corresponding to an ideal of log-quasiadjunction ${\cal A}''$
with the property described in {\it 1.} 
For any ideal $\B$ in the local ring of a singular 
point let $e_k(\B)=min_{\phi \in \B} ord_{E_k}(\pi^*(\phi))$.
Note that $\phi \in {\cal A}''$ (resp.  $\phi \in {\cal A}$)
if and only if
\begin{equation}
ord_{E_k}(\pi^*(\phi)) \ge e_k({\cal A}')
\label{order}
\end{equation} 
(resp.  $ord_{E_k}(\pi^*(\phi)) \ge e_k({\cal A})$) for all $k$ since
if $\phi \in {\cal A}''$ it certainly satisfies (\ref{order}) for each $k$ 
and vice versa if $\phi$ satisfies (\ref{order}) for all $k$ it also 
satisfies (\ref{logquasiadj}) for any $k$ and hence $\phi$ belongs to ideal 
${\cal A}''$. The same works for  $\cal A$.  
We claim that ${\cal P}({\cal A}'')$ with 
property {\it 1} is the subset of the unit cube 
which consists of solution of the
system  of inequalities in $(x_1,..,x_r)$:
 \begin{equation}
a_{k,1}{x_1}+....,+a_{k,r}{x_r}
 \ge   a_{k,1}+...+a_{k_r}-e_k({\cal A}'')-c_k-1
\label{polytope1}
\end{equation}

\noindent In order to derive ${\cal A}'' \subset {\cal A}''(j_1,..,j_r \vert
m_1,...,m_r)$ for (\ref{vector}) satisfying (\ref{polytope1}) 
note that if array $(j_1,...,j_r \vert m_1,...,m_r)$ satisfies 
(\ref{polytope1}) for all $k$ 
and if $\phi \in {\cal A}''$ i.e. we have $e_k(\phi) \ge e_k({\cal A}'')$ 
then  
we also have (\ref{logquasiadj}) for all $k$ and hence 
$\phi \in {\cal A}''(j_1,...,j_r \vert m_1,..,m_r)$. 
Vice versa, if ${\cal A}''(j_1,...,j_r \vert m_1,..,m_r)$ contains
${\cal A}''$ then for any $\phi \in {\cal A}''$ and $k$ we have 
(\ref{logquasiadj}) and hence $min_{\phi \in {\cal A}''} e_k(\phi)$
satisfies the same inequality. This proves the first part of 
the proposition.

\bigskip \noindent For the second part, let us notice that the boundary of 
the set of solutions of the system (\ref{polytope1}) is the 
set of vectors (\ref{vector}) satisfying (\ref{polytope1})
for a proper subset ${\cal S}'$ of the set of exceptional curves 
$\cal S$ and the inequalities 
\begin{equation}
a_{k,1}{x_1}+....,+a_{k,r}{x_r}
 >   a_{k,1}+...+a_{k_r}-e_k({\cal A}'')-c_k-1
\label{polytope2}
\end{equation}
\noindent for $k \in {\cal S}-
{\cal S}'$. 
For an array $(j_1,..,j_r \vert m_1,..,m_r)$, with the 
corresponding vector (\ref{vector}) in the boundary of  ${\cal P}({\cal A}'')$, let 
$E_k$ be a component with $k \in {\cal S}-{\cal S}'$. 
If $\phi \in {\cal A}''$ is a 
germ  such that $e_k({\phi})$ yields equality in (\ref{logquasiadj}) then 
the form $\omega_{\phi}$ has pole of order exactly 
one along $E_k$. Hence $ \phi \not\in {\cal A}(j_1,..,j_r \vert
m_1,..,m_r)$. Vice versa, if ${\cal A}(j_1,..,j_r \vert m_1,..,m_r)
\not= {\cal A}''(j_1,..,j_r \vert m_1,..,m_r)$ then there exist 
a form $\omega_{\phi}$ having a pole of order exactly one 
along a component $E_k$. Hence, since 
$e_k(\phi)=min \{ ord_{E_k} {\psi} \vert \psi 
\in {\cal A}''(j_1,..,j_r \vert m_1,...,m_r)\}$, the  corresponding 
vector (\ref{vector}) belongs to the boundary of the polytope 
of ${\cal A}''(j_1,..,j_r \vert m_1,..,m_r)$.
 
The vectors (\ref{vector}) corresponding to 
arrays having distinct ideals ${\cal A}$ and ${\cal A}'$
belong to faces 
which are not in the intersection of a pair of codimension one faces 
corresponding to intersecting exceptional curves (i.e. vectors (\ref{vector})
for which one has equality in (\ref{polytope1}) for  
a pair of indices such that $E_k \cap E_l \ne \emptyset$). 

\bigskip \noindent Finally, 3 follows  
since the inequalities imposed by (\ref{vector}) on $ord_{E_k} \phi$ and defining 
the ideals of quasiadjunction are the same for vectors in the interior of each face of quasiadjunction. 

\begin{proposition} Any ideal of quasiadjunciton is an ideal of 
log-quasiadjunciton (but for a different array $(j_1,..,j_r \vert m_1,..,m_r)$) 
and vice versa.
\label{comparison}
\end{proposition}

{\bf Proof.} Let ${\cal A}={\cal A}(j_1,..,j_r \vert m_1,..,m_r)$ be 
an ideal of quasiadjunciton and let $(j_1',...,j_r' \vert m_1',..,m_r')$
be an array such that the corresponding vector (\ref{vector}) belongs to 
the boundary of the set of solutions of system of inequalities 
(for any $k \in {\cal S}$):
\begin{equation}
a_{k,1}{x_1}+....,+a_{k,r}{x_r}
 \ge    a_{k,1}+...+a_{k_r}-e_k({\cal A})-c_k-1
\label{polytope3}
\end{equation}
We claim that ${\cal A}={\cal A}''(j_1,'...,j_r' \vert m_1',...,m_r')$. Indeed
if $\phi \in {\cal A}$ then $ord_{E_k}(\phi) \ge e_k({\cal A})$ together
with (\ref{polytope3}) yields 
$\phi \in {\cal A}'(j_1',..,j_r' \vert m_1',..,m_r')$.
To get opposite inclusion ${\cal A}'(j_1',..,j_r' \vert m_1',..,m_r') \subseteq
{\cal A}(j_1,...,j_r \vert m_1,..,m_r)$ notice that the vector 
(\ref{vector}) corresponding to $(j_1,...,j_r \vert m_1,....,m_r)$ 
is in the interior of the set of solutions of 
(\ref{polytope3}) since this vector cannot satisfy equality in the
system (\ref{polytope3}) because otherwise for $\phi$ such that 
$e_k(\phi)=e_k({\cal A})$ we shall have 
$a_{k,1}{{j_1+1} \over {m_1}}+....,+a_{k,r}{{j_r+1} \over {m_r}}
 =  a_{k,1}+...+a_{k_r}-e_k({\phi})-c_k-1$ contradicting 
to $\phi \in {\cal A}(j_1,...,j_r \vert m_1,..,m_r)$.
Therefore for $\phi \in 
{\cal A}''(j_1',...,j_r' \vert m_1',..,m_r')$ we shall have: 
$e_k(\phi) \ge -a_{k,1}{{j_1'+1} \over {m_1'}}-....,-a_{k,r}{{j_r'+1}
  \over {m_r'}}
 +a_{k,1}+...+a_{k_r}-c_k-1 >
-a_{k,1}{{j_1+1} \over {m_1}}-....,-a_{k,r}{{j_r+1} \over {m_r}}
 +a_{k,1}+...+a_{k_r}-c_k-1$ i.e. 
$\phi \in {\cal A}(j_1,..,j_r \vert m_1,..,m_r)$.

\bigskip 

Now let us show that any ideal of log-quasiadjunction, 
say \newline 
${\cal A}''(j_1',..,j_r' \vert m_1',...,m_r')$, is an
ideal of quasiadjunction. Let us choose array $(j_1,..,j_r \vert 
m_1,..,m_r)$ so that for corresponding vector (\ref{vector})
none of intervals 
$-a_{k,1}{{j_1'+1} \over {m_1'}}-....,-a_{k,r}{{j_r'+1}
  \over {m_r'}}
 +a_{k,1}+...+a_{k_r}-c_k-1 > x >
-a_{k,1}{{j_1+1} \over {m_1}}-....,-a_{k,r}{{j_r+1}
  \over {m_r}}
 +a_{k,1}+...+a_{k_r}-c_k-1$ contains an integer for all $k$. 
Then $\phi \in {\cal A}''(j_1',..,j_r' \vert m_1',..,m_r')$
is equivalent to $\phi \in {\cal A}(j_1,..,j_r \vert m_1,..,m_r)$.

\bigskip

The following description of the ideals of quasiadjunction 
is useful for explicit calculations of the polytopes 
introduced above and their faces.

\begin{proposition} Let $\pi: V \rightarrow {\bf C}^2$ be a 
composition of blow ups with the exceptional set $E_1 \cup ...
\cup E_k$ such that $\pi(\bigcup E_i)$ is the origin $O$. 
For a sequence of positive integers $\alpha_1,...,\alpha_k$
let $I(\alpha_1,..,\alpha_k)=\{ \phi \in {\cal O}_O \vert 
ord_{E_i} \pi^*(\phi) \ge \alpha_i, i=1,..k \}$. 

1. There are germs $\psi_i \in {\cal O}_O$ such that   
$I(\alpha_1, ...\alpha_k)$ consists of $\phi$'s such that 
the intersection index of $\phi=0$ and $\psi=0$ is  
not less than $\alpha_k$.

2. For $e_i({\cal A}'')$ in (\ref{order}) we have the identity 
${\cal A}''=I(e_1({\cal A}''),...,e_k({\cal A}''))$.
Let for such ${\cal A}''$ we have 
${\cal A}''={\cal A}''
(j_1,...,j_r \vert m_1,...,m_r)$
and the faces of 
the polytope of quasiadjunction containing corresponding
vector (\ref{vector}) are the faces corresponding to exceptional curves
$E_{k}$ where $k \in {\cal S}'$.  
Then 
\newline ${\cal A}(j_1,..,j_r \vert m_1,..,m_r)=I(...,e_k({\cal A}'')+
\epsilon_k,...)$ where $\epsilon_k=1$ for $k \in {\cal S}'$ and 
$\epsilon_k=0$ otherwise.

3. ${\cal A}'(j_1,...,j_r \vert m_1,...,m_r)$
 is the ideal generated by the ideals 
$I(e_1({\cal A}'')+\epsilon_1,...,e_k({\cal A}'')+\epsilon_k)$
corresponding all possible arrays $(\epsilon_1,...,\epsilon_k)$ 
where  $\epsilon_i=1$  if $i \in {\cal S}'$ is such 
that there exist $j \in {\cal S}'$ with $E_i \cap E_j \ne 0$ and  $\epsilon_i=0$ 
otherwise. 

\end{proposition}

\noindent {\bf Proof.} As $\psi_i$ one can take the local equation of 
the image of a transversal to $E_i$ in its generic point.
The rest follows from (\ref{order}) and the definitions.

\begin{remark} The polytopes of quasiadjunction in this paper are 
somewhat different than the local polytope of quasiadjunction of \cite{Abcov}.
The latter polytopes are defined as the equivalences classes of 
vectors (\ref{vector}) when two vectors are considered equivalent 
if and only if the corresponding ideals of quasiadjunction 
are the same. The interior of a polytope ${\cal P}({\cal A})$ is the 
union of the local polytopes of quasiadjunction from \cite{Abcov}.
Vice versa the convex polytopes in this section determine the 
local polytopes of quasiadjunction in \cite{Abcov}
\end{remark} 

\begin{remark} 
\label{multiplier}
Recall that  for a ${\bf  Q}$-divisor $D$ on a non singular
manifold $X$ its multiplier ideal  
${\cal J}(D)$ (cf. (\cite{nadel}) can be defined as follows. Let $f: Y \rightarrow X$ be an embedded 
resolution of $D$ and $f^*(D)=-E$. Then 
${\cal J}(D)=f_*({\cal O}_Y(K_Y-f^*(K_X)-\lfloor E \rfloor))$ where $\lfloor E \rfloor$
is round-down of a $\bf Q$-divisor. In this terminology one can define the ideals of 
quasiadjunction as follows. For an array $(\gamma_1,..,\gamma_r), (\gamma_i \in {\bf Q})$ 
let $D_{\gamma_1,..,\gamma_r}$ be given by equation $f_1^{\gamma_1} \cdot \cdot \cdot 
f_r^{\gamma_r}$. Then ${\cal J}(D_{\gamma_1,..,\gamma_r})={\cal A}(j_1,..,j_r \vert m_1,...,
m_r)$ where $\gamma_i=1-{{j_i+1} \over {m_i}}$ for $i=1,..,r$. 
This follows immediately from (\ref{quasiadj}).
\end{remark}
\subsection{Mixed Hodge Structure of the cohomology of links of 
singularities and on the local cohomology.}
\label{mixedhodge}

Here we shall summarize several well known facts used in the
next section. 
Cohomology of the link $L$ of an isolated singularity $x$ of a 
complex space $X$ (${\rm dim}X=n$) can be given a Mixed Hodge structure for example 
using canonical identification $H^k(L)=
H^*_{\{x\} }(X)$ with the local cohomology 
and using the construction of mixed Hodge structure on the latter due
to Steenbrink \cite{arcata}. 
The Hodge numbers: $h^{kpq}(L)=
{\rm dim} Gr^p_FGr^W_{p+q}H^k(L)$ have the following 
symmetry properties (cf. \cite{hamm},
 \cite{arcata}):
\begin{equation}
h^{kpq}=h^{2n-k-1,n-p,n-q}
\label{symmetry}
\end{equation} 
If $E$ is the exceptional divisor for a resolution then
for $k <n$ one has (cf. \cite{hamm}) 
$$
h^{kpq}(L)=h^{kpq}(E) \ \ \ if \ p+q <k$$
\begin{equation}
h^{kpq}(L)=h^{kpq}(E) - h^{2n-k,n-p,n-q}(E)
\ \ \ if \ p+q =k
\label{reduction1}
\end{equation} 
$$h^{kpq}(L)=0 \ \ \ if \ p+q >k$$
The mixed Hodge structure on cohomology of a link is related to the 
mixed Hodge structure on vanishing cohomology via 
the exact sequence (corresponding to the exact sequence of a pair):
\begin{equation}
0 \rightarrow H^{n-1}(L) \rightarrow H^n_c(B) \rightarrow H^n(B)
\rightarrow H^n(L) \rightarrow 0
\label{reduction2}
\end{equation}
which is an exact sequence of mixed Hodge structures (cf. \cite{oslo} and (2.3) in \cite{arcata}).
\par Steenbrink also put Mixed Hodge structure of the local 
cohomology $H^*_E(\tilde X)$ (\cite{arcata}) where $\tilde X$ is a 
resolution of $X$. 
In the case $dim_{\bf C} \tilde X=2$ 
we have 
\begin{equation}
H^*_E(\tilde X)=Hom(H^{4-*}(E),{\bf Q}(-2))
\label{duality}
\end{equation}
where ${\bf Q}(-2)$ is Tate Hodge of type $(2,2)$. 
Since the Hodge and weight filtrations on 
$H^1(E)$ have form: 
$$H^1(E)= W_1 \supset W_0 \supset 0, \ \ \ H^1(E)=F^0 \supset F^1 \supset
F^2=0$$
we have on $H^3_E(\tilde X)$
 
$$H^3_E(\tilde X)=W_4 \supset W_3 \supset W_2=0,H^3_E(\tilde X)=F^1 \supset F^2 \supset
F^3=0$$

Moreover 
\begin{equation}
F^1H^1(L)=F^1H^1(E)=F^2H^3_E(\tilde X)
\label{hodgefiltration}
\end{equation}

One can use the following complex for description of this 
mixed Hodge structure: 

\begin{equation}
0 \rightarrow A^2_E(\tilde X) \rightarrow A^3_E(\tilde X) 
 \rightarrow 0
\label{complex}
\end{equation}

where $$A^2_E(\tilde X)=
\Omega^1_{\tilde X}(log \ E)/\Omega^1_{\tilde X},
A^3_E(\tilde X)=
\Omega^2_{\tilde X}(log \ E)/\Omega^2_{\tilde X}$$

with filtrations given by 
$$F^2A^p_E(\tilde X)=0 \ for \ p<3,
F^2A^p_E(\tilde X)=A^p_E( \tilde X) \ for \ p \ge 3$$

$$W_3A^3_E(\tilde X)=W_1\Omega^2_{\tilde X}(log \ E)/\Omega^2_{\tilde X}$$

Since $H^3(E)=0$, the relations (\ref{reduction1}) and (\ref{duality}) yield 
that the complex (\ref{complex}) completely determines 
$h^{1pq}$ (and hence all Hodge numbers $h^{kpq}$ by (\ref{symmetry})).

\section{Characteristic varieties and polytopes of quasiadjunction.}

\subsection{Main Theorem}
We shall view the unit cube 
$\cal U$, considered in the last section and 
containing the polytopes of quasiadjunction 
as the fundamental domain for the Galois group 
$H^1(S^3-L,{\bf Z})$
of the universal 
abelian cover $H^1(S^3-L,{\bf R})$
 of the group $H^1(S^3-L,{\bf R}/{\bf Z})$
of the unitary characters of 
$H_1(S^3-L,{\bf Z})$
(i.e. the maximal compact subgroup 
of $Char (H_1(S^3-L,{\bf Z}))=H^1(S^3-L,{\bf C}^*)$).
$exp: {\cal U} \rightarrow 
Char (H_1(S^3-L,{\bf Z}))$ will denote the restriction of 
$H^1(S^3-L,{\bf R}) \rightarrow H^1(S^3-L,{\bf R}/{\bf Z})$
on $\cal U$.
\par For any sub-link $\tilde L$ of $L$ i.e. a link formed  by 
components of $L$ we have surjection $\pi_1(S^3-L) \rightarrow 
\pi_1(S^3-\tilde L)$ induced by inclusion. Hence $Char H_1(S^3-\tilde L,{\bf Z})$
is a sub-torus of $Char H_1(S^3-L,{\bf Z}))$ (in coordinates in the latter 
torus corresponding to the components of $L$ it is given by equations 
of the form $t_{\alpha}=1$ where subscripts correspond to components of $L$
absent in $\tilde L$). Moreover, since the homology of the universal abelian cover 
$H_1(\widetilde {S^3-L})$ surjects onto
$H_1(\widetilde {S^3-\tilde L})$ it follows that $V_i(S^3-\tilde L)$ belongs to a 
component of $V_i(S^3-L)$ (cf. 1.2.1 in \cite{Abcov}). 
We shall call a character of $\pi_1(S^3-L)$ (or a connected component of $V_i(S^3-L)$) 
{\it essential} if it does not belong to a subtorus $Char H^1(S^3 -\tilde L)$ 
for any sublink $\tilde L$ of $L$.
\par Let $L_{m_1,..,m_r}$ be the link of singularity
(\ref{completeintersection}) or equivalently the cover 
of $S^3$ branched over link $L$ and having a quotient 
$H_{m_1,..,m_r}={\bf Z}/m_1{\bf Z} \oplus .... \oplus 
{\bf Z}/m_1{\bf Z}$ of  $H_1(S^3-L,{\bf Z})$ 
as its Galois group. We shall view $Char H_{m_1,..,m_r}$
as a subgroup of $Char H_1(S^3-L,{\bf Z})$. The group 
$H_{m_1,..,m_r}$ acting on 
$H^1(L_{m_1,..,m_r})$ preserves both Hodge and weight filtrations.
\noindent \begin{theorem} An essential  character $\chi \in 
Char (H_1(S^3-L,{\bf Z})) $ is a character of 
the representation of $H_{m_1,..,m_r}$ acting on 
$F^1(H^1(L_{m_1,..,m_r}))$ if and only if it factors through
the Galois group $H_{m_1,..,m_r}$ and belongs to the image
of a face of quasiadjunction under the exponential map.
\par The multiplicity of $\chi$ in this representation of the
Galois group is equal to $dim {\A}_{\Sigma}''/{\cal A}_{\Sigma}$
where $  {\cal A}_{\Sigma}''$ (resp. ${\cal A}_{\Sigma}$)
is the ideal of log-quasiadjunction (resp. ideal of quasiadjunction)
corresponding to a vector (\ref{vector}) belonging to 
the face of quasiadjunction $\Sigma$ 
\par A character $\chi$ is a character of
the representation of the Galois group of the cover on
$W_0(H^1(L_{m_1,..,m_r}))$ if and only if it factors through
the Galois group $H_{m_1,..,m_r}$ and it belongs to the image
under the exponential map of a weight one face of quasiadjunction.

\label{main}
\end{theorem}

\proof {\it 1. log-2-forms on $V_{m_1,...,m_r}-p$}.
Let
$\tilde p: \tilde V_{m_1,...,m_r} \rightarrow V_{m_1,..,m_r}$
be a resolution such that the exceptional locus is a  divisor
$\tilde E=\bigcup \tilde E_i$ on
$\tilde V_{m_1,...,m_r}$ with normal crossings (e.g. (\ref{resolution})). 
The group
$H_{m_1,...,m_r}$ acts on both sheaves: 
\newline $\Omega^2_{\tilde V_{m_1,...,m_r}}
(log \ E)$ and $\Omega^2_{\tilde V_{m_1,..,m_r}}$.

We are going to identify the eigenspace of 
$\Omega^2_{\tilde V_{m_1,...,m_r}}
(log \ E) / \Omega^2_{\tilde V_{m_1,..,m_r}}$
corresponding to the character $\chi$, which 
is the exponent of a vector (\ref{vector}) belonging to a face of
quasiadjunction $\Sigma$, with the quotient of ideals
${\cal A}_{\Sigma}''/{\cal A}_{\Sigma}$.
\par First notice that any 2-form ${z_1^{j_1} \cdot \cdot \cdot z_r^{j_r}
\phi(x,y) dx \wedge dy} \over {z_1^{m_1-1} \cdot \cdot \cdot
z_r^{m_r-1}}$ is holomorphic on $\tilde V_{m_1,..,m_r}-E=
V_{m_1,...,m_r}-p$ since 
it is a residue of a holomorphic $(r+2)$-form on ${\bf C}^{r+2}-V_{m_1,...,m_r}$.
It is an eigenform corresponding to the character $\chi$ such that $\chi(\gamma_i)=
exp(2 \pi \sqrt {-1} { {m_i-j_i-1} \over {m_i}})$.
\par Vice versa any 2-form, holomorphic on $V_{m_1,...,m_r}-p$, is a residue of
$r+2$-form given by ${\phi(z_1,..,z_r,x,y) dz_1 \wedge ... \wedge
dz_{r} \wedge dx \wedge dy} \over {(z_1^{m_1}-f_1(x,y)) \cdot \cdot
\cdot (z_r^{m_r}-f_r(x,y))}$ where $\phi$ is a polynomial.
Decomposition of $\phi$ into a sum of monomials corresponds to
the decomposition of the form into sum over characters.
\par Secondly a form ${z_1^{j_1} \cdot \cdot \cdot z_r^{j_r} \phi
dx \wedge dy} \over {z_1^{m_1-1} \cdot \cdot \cdot z_r^{m_r-1}}$
is log-2-form (resp. holomorphic 2-form) 
if and only if $\phi(x,y)$ is in the ideal of 
log-quasiadjunction (resp. quasiadjunction) corresponding to 
$({{j_1} \over {m_1}},....,{{j_r} \over
{m_r}})$.
These ideals do not coincide if and only if  (\ref{vector})
belongs to a face of quasiadjunction.

\par {\it 2. Hodge and weight filtration on cohomology 
of link.} 
Now we want to identify  
$\Omega^2_{\tilde V_{m_1,...,m_r}}(log \ E)/\Omega^2_{\tilde V_{m_1,...,m_r}}
$ with $F^1 H^1(L_{m_1,...,m_r})$.
(\ref{hodgefiltration}) yields that the latter is isomorphic to
$F^2H^3_E(\tilde V_{m_1,...,m_r})$. Using  description of the
Hodge filtration on $H^3_E(\tilde V_{m_1,...,m_r})$ from  (\ref{mixedhodge})
it can be identified with the hypercohomology of the complex
$0 \rightarrow A^3_E(\tilde V_{m_1,...,m_r}) \rightarrow 0$ i.e. with 
$H^0(\Omega_{\tilde V_{m_1,...,m_r}}^2(log \ E)/\Omega_{\tilde V_{m_1,...,m_r}}
^2)$.
Since by Grauert-Riemenschnider theorem $H^1(\Omega^2_{\tilde V_{m_1,...,m_r}})=0$ we see that 
the latter space is isomorphic to 
$H^0(\Omega_{\tilde V_{m_1,...,m_r}}^2(log \ E))/H^0(\Omega_{\tilde V_{m_1,...,m_r}}^2)$ 
and the claim follows.

\par {\it 3. Conclusion of the proof.} Similarly to 
the above, it follows that 
a character $\chi$ is a character of the representation of $H_{m_1,..,m_r}$ on 
$W_0H^1(L_{m_1,..,m_r})$ if it is a character of this  group acting on    
$W_2/W_1(\Omega^2_{\tilde V_{m_1,..,m_r}}(log \tilde E)/
\Omega^2_{\tilde V_{m_1,..,m_r}})$ (cf. (\ref{mixedhodge})). 
In other words
$\chi=(exp({{2 \pi i j_1} \over m_1},..., {{2 \pi i j_r} \over m_r})$
where ${\cal A}'(j_1,...,j_r \vert m_1,...,m_r) \ne 
{\cal A}''(j_1,...,j_r \vert m_1,...,m_r)$ i.e. 
$({{j_1} \over {m_1}},....,{{j_r} \over
{m_r}})$ belongs to a face of weight one log quasiadjunction.
Moreover, the dimension of the $\chi$-eigenspace of the action of the 
Galois group of the cover on $W_0H^1(L_{m_1,..,m_r})$ is equal to 
$dim {\cal A}''_{\Sigma}/{\cal A}'_{\Sigma}$ where $\Sigma$ is 
the face of weight one log-quasiadjunction to which the  
$\chi$ belongs.

\subsection{Essential components of characteristic varieties.}

Theorem \ref{main} allows to describe essential components of 
characteristic varieties. Indeed, each component of 
$V_i(S^3-L)$ is a torus translated by a point of a finite order 
in $Char H_1(S^3-L,{\bf Z})$ (cf. \ref{translatedtori}) and each 
such sub-torus is Zariski closure of the set of points of finite
order in it. It follows from \cite{sakuma} that a essential character 
of finite order belongs to $V_i(S^3-L)$ is and only if it is a character of 
$H_{m_1,..,m_r}$ on $H_1(L_{m_1,..,m_r})$ for some array $(m_1,..,m_r)$. 
A character $\chi$ appears as either a character on $W_0H^1(L_{m_1,..,m_r})$, 
in which case  according to \ref{main} it is exponent of a vector in a
face of quasiadjunction, or $\chi$ is a character on $W_1/W_0H^1(L_{m_1,..,m_r})$,
in which case either $\chi$ or $\bar \chi$ is a character of the 
Galois group acting on $H^{1,0}(W_1/W_0)$. 
In each of the cases the multiplicity of the character is 
${\cal A}''_{\Sigma}/{\cal A}_{\Sigma}$ 
where $\Sigma$ is the face of quasiadjunction to exponent 
of which $\chi$ belongs. If ${\cal L}_{\chi}$ will denote
the local system on $S^3-L$ corresponding to the character 
$\chi$ then this multiplicity is equal to 
to $dim H^1({\cal L}_{\chi})$,
as follows from arguments 
in \cite{sakuma}. On the other $dim H^1({\cal L}_{\chi})$
is the depth of the characteristic variety to which $\chi$ belongs.
We obtain therefore:

\begin{proposition} For $\xi=(x_1,..,x_r) \in {\cal U}$ let $\bar \xi=(1-x_1,....1-x_r)$.
For a face of quasiadjunction $\Sigma$ let $\tilde \Sigma=\Sigma$ 
or $\tilde \Sigma=\Sigma \cup \bar \Sigma$ depending on whether  
${\cal A}_{\Sigma}'/{\cal A}_{\Sigma}=0$ or not. Then Zariski closure of 
$exp(\tilde \Sigma)$ is a component of characteristic variety of depth 
$dim {\cal A}_{\Sigma}'/{\cal A}_{\Sigma}+dim {\cal A}_{\bar \Sigma}'/{\cal A}_{\bar \Sigma}
+dim {\cal A}_{\Sigma}''/{\cal A}_{\Sigma}'$. 
Vice versa, any essential component 
of $V_i(S^3-L)$ has such form.

\end{proposition}

\begin{subsection}{Quasiadjunction and higher Alexander polynomials.}
\bigskip
We shall show that in the case $r=1$ information from the 
 faces and ideals of quasiadjunction determines  
{\it all} Hodge numbers $h^{p,q}_{\zeta}$ of the Milnor fiber of 
$f=0$ considered in introduction 
and in particular {\it all} Alexander polynomials $\Delta_i$  
(cf. Introduction).
The idea to relate the constants of quasiadjunction to the 
Hodge theory appeared first in \cite{loser} where the case of 
the first Alexander polynomial $\Delta_1$ 
was studied.  

For $r=1$ each polytope of quasiadjunction is a segment $[\kappa,1]$.
$\kappa$ is called (cf. \cite{arcata1}) the constant of
quasiadjunction (resp. log-quasiadjunction and weight one quasiadjunction).
  
\begin{proposition}
\par 1. $\zeta$ is an  eigenvalues of the monodromy acting on the 
Milnor fiber of $f$ if and only if 
$\zeta$ or $\bar \zeta$ is equal to $exp 2 \pi \kappa_l$ for some 
constant of quasiadjunction $\kappa_l$.
\par 2. If $\zeta=exp 2 \pi \kappa_l$ and $\kappa_l$ is 
a weight one constant quasiadjunction then 
$h^{0,0}_{\zeta}=
dim {\cal A}_{\kappa_l}''/{\cal A}_{\kappa_{l}}'$ and
 $h^{1,0}_{\zeta}=dim {\cal A}_{\kappa_l}'/{\cal A}_{\kappa_{l}}$.
In particular $h^{0,0}_{\zeta}=0$ unless corresponding $\kappa_l$ 
is a weight one constant of quasiadjunction.
\end{proposition}

\noindent \proof $\zeta \ne 1$ such that $\zeta^m=1$ is an eigenvalue of 
the monodromy of the Milnor fiber of $f$ if an only if it is an eigenvalue of 
the Galois group of acting on $H^1(L_m)=H^1_{\{SingV_m\}}(V_m)$
where $V_m$ is given by $z^m=f$ (cf. (\ref{completeintersection})).
On the other hand, the exact sequence 
$$0 \rightarrow H^1_{\{SingV_m\}}(V_m) \rightarrow H^2_c(M(V_m)) 
\buildrel j \over \rightarrow H^2(M(V_m)) \rightarrow $$
(cf. (\ref{reduction2})) in which $M(V_m)$ is the Milnor fiber of 
$V_m$ shows that the Hodge structure on $H^1(L_m)$ is isomorphic to the 
Hodge structure induced from $H^2_c(B)$ on $Ker(T_c-I)$.
Indeed $ker \ j=ker \ T_c-I$ since 
$T_c-I =Var \circ j$ (cf. \cite{arcata}, (2.4); recall also that  
$Var: H^2(B) \rightarrow H^2_c(B)$ is an isomorphism as a consequence for 
example of a well known relation between the Seifert form 
and variation operator). 

The result follows since the Hodge structure on $H^2_c(M(V_m))$ is determined by the Hodge 
structure on the Milnor fiber  of $f=0$ via Thom-Sebastiani type 
theorem (cf. \cite{oslo}).

\end{subsection}

\section{Properties and applications of polytopes of quasiadjunction}
\label{properties}

\subsection{Semicontinuity.}
\label{semicontinuity}
\bigskip  
\begin{theorem} 1.Let $C_t$ be a family of plane curves singularities
in a ball $B$ with $r$ branches such that the limit curve has $r$ branches as well.
Then the number of components of $V_1(C_t)$ does not exceed the number
of components of $V_1(C_0)$.  
\par 2. Total volume of all codimension one faces of quasiadjunction 
is semicontinuous  provided that the volume of each face is calculated 
with respect to the measure on the hyperplane containing this face 
in which the measure of the simplex containing no integral points 
is equal to $ 1 \over {(r-1)!}$.
\end{theorem}

\proof  First notice that the intersection form on $H_2$ of the smoothing of 
complete intersection surface singularity which is an abelian
cover of $B$ branched over $C_t$ (i.e the singularity (\ref{completeintersection}))
embeds into the intersection form of 
the singularity which is the abelian cover 
of the same type branched over $C_0$. On the other hand,  
$b_1(L_{m_1,...m_r})$ is the dimension of the radical of this 
intersection form. In particular we have $b_1(L_{n,...n}(C_t)) \le b_1(L_{n,...n}(C_0))$.
On the other hand $b_1(L(n,...n))=kn^{r-1}+\alpha \cdot n^{r-1}+...$ 
for almost all $n$ as follows from \cite{sakuma}
where $k$ is the number of components in $V_1(C)$ since 
the number of points of order $n$ on a torus of dimension $l$ translated by a point
of finite order is $n^l$ for almost all $n$.
Hence $k=\lim_{n \rightarrow \infty} {b_1(n,..n) \over {n^{r-1}}}$ which 
yields the first part.
\par Similarly,  asymptotically the number of the points in the
lattice $({k \over n} .... {k \over n}) \subset {\cal U}$ which belong 
to codimension one faces of quasiadjunction is the total volume of the faces.
On the other hand this number of the points is $dim F^1 H^2_1(M)$ i.e. the 
dimension of the Hodge filtration $F^1$ on the subspace of the 
cohomology of the Milnor fiber consisting of monodromy invariants.
Similarly to \cite{semicontinuity} for the Milnor 
fibers $M_t$ and $M_0$ of the smoothings of abelian covers 
of $B$ having type $(n,...n)$ and branched over $C_t$ and $C_0$ respectively we have:
$dim F^1 H^2_1(M_t)-dim F^1H^2_1(M_0)=dim F^1 {\bf H}_1(R\Phi)$ 
which yields the second part.

\subsection{Log canonical divisors} 
\label{thresh}

Recall (\cite{kollar})
that a pair $(X,D)$ where $X$ is normal and $D$ is 
a $\bf R$-divisor such that $K_X+D$ is $\bf R$-Cartier 
is called log-canonical at $x \in X$
if for any birational morphism
$f: Y \rightarrow X$, with $Y$ normal, in the decomposition 
\begin{equation} K_{Y}=f^*(K_X+D)+\sum_E a(E,X,D)E
\label{discrepancy}
\end{equation} 
for each irreducible $E$ having center at $x$ one has $a(E,X,D) \ge -1$. This 
coefficient is called {\it discrepancy} of divisor $D$ 
on $X$ along $E$.

\begin{proposition} The local ring ${\cal O}_O$  
of a singularity $f_1 \cdot \cdot \cdot f_r=0$ at the origin $O$
of ${\bf C}^2$ considered as the ideal in itself 
is an ideal of log-quasiadjunction.
Let $\cal P$ be the corresponding polytope of log-quasiadjunction.
Let $D_i$ be the divisor in ${\bf C}^2$ with the local equation $f_i=0$ 
near the origin. 
\newline \noindent Then for $\{(\gamma_1,...,\gamma_r ) \} \in  {\bf R}^r$
the divisor $\gamma_1 D_1+...+\gamma_r D_r$ is log-canonical 
at $(0,0) \in {\bf C}^2$ if and only if $(\gamma_1+1,..,\gamma_r+1)$
belongs to the polytope $\cal P$.
\label{treshold}
\end{proposition}

\noindent \proof Let us consider the polytope given 
by inequalities (\ref{polytope3}) in which one puts 
$e_k({\cal A}'')=0$ i.e. 
 \begin{equation}
a_{k,1}{x_1}+....,+a_{k,r}{x_r}
 \ge   a_{k,1}+...+a_{k_r}-c_k-1
\label{polytope4}
\end{equation}
 Let $(j_1,..,j_r \vert m_1,..,m_r)$
be such that the corresponding vector (\ref{vector})
belongs to the boundary of this polytope. 
Then $1 \in {\cal A}''(j_1,..,j_r \vert m_1,..,m_r)$ and hence
${\cal A}''(j_1,..,j_r \vert m_1,..,m_r)$
is the local ring of the origin. 
\par If $\pi: V \rightarrow {\bf C}^2$ is an embedded resolution
$\pi^*(f_1^{\gamma_1} \cdot \cdot \cdot f_r^{\gamma_r} dx \wedge 
dy)$ has as the order of vanishing along $E_k$:
  $$a_{k,1}\gamma_1+...+a_{k,r} \gamma_r +c_k$$
i.e. the discrepancy along each $E_k$ is not less than $-1$
if and only if $(\gamma_1+1,...,\gamma_r+1)$ satisfies
(\ref{polytope4}). 

\begin{remark} If $r=1$ the polytope of quasiadjunction 
$\cal P$ described above has as its face (i.e. the end)
the constant of quasiadjunction (cf. \cite{arcata1} for a definition)
$\kappa_1$ and one obtains
that the log-canonical threshold $c_0(f)$ of $f=0$ satisfies:
$c_0(f)+1=\kappa_1$. This result is due to J.Kollar (cf. (\cite{kollar})
Prop. 9.8).
\end{remark}

\begin{remark}{{\rm Mixed Hodge structure on cohomology of 
universal abelian covers.}} It seems at the moment very little 
is known about the Hodge theory on cohomology of complex manifolds
which are not finite CW-complexes. Theorem \ref{main} 
gives some hints what one may and may not expect. 
It yields a decomposition of the torus of unitary 
characters into a union of connected subsets 
(faces of quasiadjunction or their conjugates)
in which the characters 
appearing on certain Hodge and weight type on finite 
level are dense. These subsets are not algebraic subvarieties. 
In particular one cannot expect filtrations on the cohomology 
of infinite abelian covers with Galois-invariant Hodge and 
weight filtrations. 
Results of \cite{Abcov} yield  
decomposition of the torus of 
unitary characters of the fundamental group of the complement
to an affine algebraic curve similar to the one given by theorem 
\ref{main} and reflecting the Hodge theory on $H^1$ 
of the infinite abelian cover of this complement.

\end{remark}

\begin{section}{Examples}

\bigskip 

{\it 1.Links of Ordinary singularities.} Let us consider the 
link of the singular point of $L_1 \cdot \cdot \cdot L_r=0$ where $L_i$ are distinct 
linear forms on ${\bf C}^2$ (the Hopf link with $r$ components).
The faces of 
quasiadjunction are 
\begin{equation}
\label{ordinaryface} \Sigma_l: \ \ \ x_1+...+x_r=l, \ \ \ (l=1,...,r-2)
\end{equation}
This follows directly since a 
resolution is achieved by single blow up and the multiplicity of  
 $L_i$ on the exceptional curve is equal to 1 (i.e. the coefficients $a_{k,i}$
in (\ref{logquasiadj}) are all equal to $1$).
Moreover, if $\cal M$ is the maximal ideal in the local ring of the origin,
then the ideal of 
quasiadjunction ${\cal A}_{\Sigma_l}$ corresponding to the face (\ref{ordinaryface})
is ${\cal M}^{r-1-l}$ and ${\cal A}_{\Sigma_l}'={\cal A}_{\Sigma_l}''={\cal M}^{r-2-l}$.
We have $dim {\cal M}^{l-1}/{\cal M}^{l}=l$. Hence, if $\chi$  
is a character of the fundamental group having multi-order 
$(m_1,..m_r)$ and 
\begin{equation} \chi(\gamma_i)=
exp (2 \pi \sqrt {-1} x_i) \ \  {\rm where} \ \Sigma x_i=l 
\label{ordinarysing}
\end{equation} 
then the dimension of the eigenspace 
$H^{1,0}(L_{m_1,..m_r})_{\chi}$ is $r-1-l$. 
The conjugate character
$\bar {\chi}$ has the same multiplicity $r-1-l$ on $H^{0,1}(L_{m_1,..m_r})$. 
Therefore the character in (\ref{ordinarysing}) has on $H^{0,1}(L_{m_1,..m_r})$ 
multiplicity $l-1$ and its multiplicity  on $H^1(L_{m_1,..m_r})$ is equal to $r-2$. 
In particular, the characteristic variety has only one essential component $V_{r-2}$. This 
is well known from the calculations using the Fox calculus (cf. \cite{topandappl}  for $r=3$).

\bigskip{\it 2.Singularity $(x^2+y^5)(y^2+x^5)$} (cf. also \cite{volodia}).
A standard sequence of blow ups 
lead to the resolution $V$ for this singularity pictured on figure 1
($E_1,E_2$ are the exceptional curves in the last blow up and 
$B_1,B_2$ are the branches of the singularity).
\begin{figure}
 \centering
 \includegraphics[height=4cm]{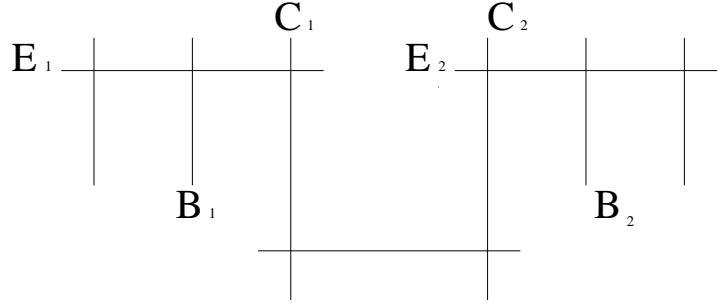}
 \caption{Resolution of $(x^2+y^5)(x^5+y^2)$.}
 \label{graph}
\end{figure}
\par We have $e_{E_1}(x)=5,e_{E_1}(y)=2,e_{E_2}(x)=2,e_{E_2}(y)=5$ and hence 
$e_{E_1}(x^2+y^5)=10, e_{E_1}(y^2+x^5)=4,e_{E_2}(x^2+y^5)=4,e_{E_2}(y^2+x^5)=10$.
Moreover $e_{E_1}(dx \wedge dy)=e_{E_2}(dx \wedge dy)=6$ (the orders of all functions 
and forms on ${\bf C}^2$ meant to be 
calculated on
the chosen  resolution $V$). As in (\ref{completeintersection}),
 $V_{m_1,m_2}$ is the abelian 
branched cover: $z_1^{m_1}=x^2+y^5, z_2^{m_2}=x^5+y^2$. 
A form $\omega_{\phi}$ given 
by (\ref{form}) admits a 
holomorphic extension over the preimages of the curves $E_1$ and $E_2$
on the normalization 
$\widetilde V_{m_1,m_2}$ of $V \times_{{\bf C}^2} V_{m_1,m_2}$
if and only if 
\begin{equation}
 10({{j_1+1} \over {m_1}}-1)+4({{j_2+1} \over m_2}-1)+e_{E_1}(\phi)+7 \ge 0
\label{exampl}
\end{equation}
$$4({{j_1+1} \over {m_1}}-1)+10({{j_2+1} \over m_2}-1)+e_{E_2}(\phi)+7
\ge  0$$
(cf. (\ref{logquasiadj})). 
Similar inequalities should be written for other exceptional curves but 
calculation shows that the inequalities expressing the condition that 
$\omega_{\phi}$ extends over 
the preimages of remaining exceptional curves in figure 1  follow from (\ref{exampl}). 
\par The quotient ${\cal A}^{\prime \prime}(j_1,j_2 \vert m_1,m_2) /{\cal A} 
(j_1,j_2 \vert m_1,m_2) \ne 0$ if and only if there exist a 
function $\phi$ for which (\ref{exampl}) is satisfied and 
such that a left hand side in at least one of inequalities 
(\ref{exampl}) is zero.
 Since the possibilities for $e_{E_i}(\phi)$ are $0,2,4,5,6,7,...$ the 
faces of quasiadjunction are subsets of the lines: 
$$ 10x_1+4x_2=7,5,3,2,1$$
$$ 4x_1+10x_2=7,5,3,2,1$$
The value $e_{E_1}(\phi)=5$ (i.e. $\phi=x$) for $i=1$ 
  does not yield a face of quasiadjunction since
it implies that $e_{E_2}(\phi)=2$ and no pair 
$x_1={{j_1+1} \over m_1},x_2={{j_2+1}\over
m_2}$ in unit square satisfies  inequalities (\ref{exampl})
with such $e_{E_i}(\phi)$. Similarly, $e_{E_2}(\phi)=5$ also does not 
yield a face of quasiadjunction.
\par  Next let us consider pairs
$({{j_1+1} \over m_1},{{j_2+1} \over m_2})$ satisfying $10x_1+4x_2=7$. 
It follows from the first inequality in  (\ref{exampl}) that $\phi$ is a non zero in  
\newline ${\cal A}^{\prime \prime} (j_1,j_2 \vert m_1,m_2)/{\cal A}(j_1,j_2 \vert m_1,m_2)$ only if 
$e_{E_1}(\phi)=0$ and hence the lowest order term of $\phi$ is a constant (i.e. this is the only
case when $\omega_{\phi}$ has pole of order 1 along $E_1$).  
Hence the second inequality in  (\ref{exampl})
 is $4x_1+10x_2 \ge 7$ and hence 
only ``the half'' of the segment $10x_1+4x_2=7$ in the unit square is
the face of quasiadjunction. Similarly ``the half''
of the segment $4x_1+10x_2=7$ is 
also the face of quasiadjunction.
\par On the other hand for a pair $({{j_1+1} \over m_1},{{j_2+1} \over m_2})$
on the segment $10x_1+4x_2=5$, 
assuming that $\omega_{\phi}$ has pole of order one along $E_1$,
 the first inequality (\ref{exampl}) yields that the lowest order term of $\phi$ 
is $ay, a \in {\bf C}^*$
and the second one in (\ref{exampl}) is satisfied for any pair 
$({{j_1+1} \over m_1},{{j_1+1} \over m_1})$ on 
$10x_1+4x_2=5$ in the unit square. Therefore the segement $10x_1+4x_2=5$ is 
a face of quasiadjunction. Simliar calculations show that
we obtain the diagram 
of faces quasiadjunction given on figure 2.
 \begin{figure}
 \centering
 \includegraphics[height=6cm]{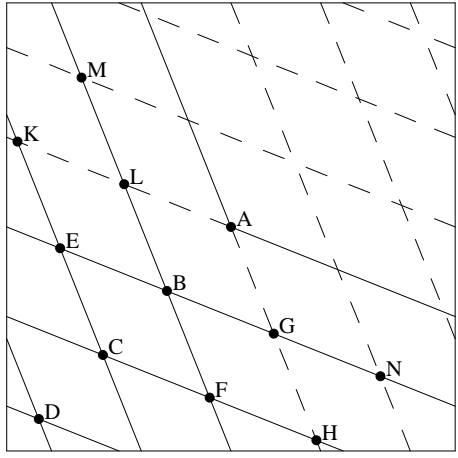}
 \caption{Faces of quasiadjunction for $(x^2+y^5)(x^5+y^2)$.}
 \label{graph}
\end{figure}
Moreover the points $B,C,D,E,F$ are the only ones for which one has 
$dim {\cal A}^{\prime  \prime}/{\cal A}=2$ (for the remaining 
points on the faces of quasiadjunction this dimension is 1).
\par  Exponential map takes the union 
of the faces of quasiadjunction and their conjugates into the union of 
of translated subgroups $t_1^2t_2^5+1=0$ and $t_1^5t_2^2+1=0$
(and coincides with the union of translates of the maximal compact subgroups
of the latter). $V_2$ consists of exponents of the points 
of following types:

\par \noindent a) points where 
$dim {\cal A}^{\prime  \prime}/{\cal A}=2$ (corresponding 
characters appears on holomorphic part); these are the points $B,C,D,E,F$.
\par \noindent b) conjugates of the points in a) (corresponding 
characters appear on anti-holomorphic part).
\par \noindent c) points belonging to faces of quasiadjunction and having 
conjugates on faces of quasiadjunction as well; 
corresponding characters appear with multiplicty
one on both the holomorphic and anti-holomorphic parts;
these are points $H,G,K,L,M,N$ ($N$ is the conjugate of $M$.) 
\par \noindent d) conjugates of the points in c)

\par \noindent We obtain that $V_2$ consists of all (twenty) points of intersection of two translates
$t_1^2t_2^5+1=0$ and $t_1^5t_2^2+1=0$ except the point $(-1,-1)$. 

The eigenspace in $H^1(L)$ corresponding to 
the character $exp (2 \pi i A)=(-1,-1)$ 
consists of weight zero classes in $H^1(L)_{(-1,-1)}$  since 
in this case ${\cal A}^{\prime \prime}
/{\cal A}^{\prime}$ is generated by $1$ and $\omega_1$ has poles of order 
1 along $E_1,E_2$ and $C_1,C_2$ i.e. weight of $\omega_1$ is two.

\end{section}

\end{document}